\documentclass[dvips,a4paper,10pt]{article}
\usepackage{epsfig,graphicx,verbatim,amsmath,amsfonts,amssymb,alltt,algorithm,algorithmic,multirow}
\usepackage[margin=1in,nohead]{geometry}

\begin{document}
\newcommand{\raf}[1]{(\ref{#1})}
\newcommand{\sch}{\choose }
\newcommand{\ul}[1]{\underline{#1}}
\newcommand{\comp}{\hbox{$<\kern -3pt >$}}
\newcommand{\ncomp}
		{\;\hbox{\hbox{/}\kern -9.5pt \hbox{$<\kern -3pt >$}}}
\newcommand{\mpt}{\partial}

\newcommand{\meet}
	       {\hbox{$\wedge \kern -5.75pt \raise 1.5pt \hbox{$.$}\,$}}
\newcommand{\Meet}
	     {\hbox{$\bigwedge \kern -8pt \raise 0.75pt \hbox{$.$}\:$}}
\newcommand{\pplus}
	    {\hbox{$+ \kern -6pt \raise 9.05pt \hbox{$.$}\:$}}
\newcommand{\ld}
	       {\hbox{$< \kern -6pt \raise 2pt \hbox{$.$}\,$}}

\newcommand{\bve}{\bigwedge}
\newcommand{\fml}{((\mu-L)/\gs)}
\newcommand{\fmgl}{((\mu-L)/(\ga\mu))}

\newcommand{\ctnt}{\hbox{$\,\hat{\ }\,$}}

\newcommand{\ct}{\centerline}
\newcommand{\sat}{\models}
\newcommand{\prvs}{\vdash}
\newcommand{\frc}{\hbox{$\parallel \kern -5.7pt \hbox{$-$}$}}
\renewcommand{\iff}{\leftrightarrow}
\newcommand{\ra}{\rightarrow}

\newcommand{\lharp}{\leftharpoonup}
\newcommand{\rharp}{\rightharpoonup}

\newcommand{\nfrc}{\not \kern -5pt \frc}
\newcommand{\rest}{\vbox{\hbox{$\:\kern -2pt\mathbin{\vert\kern-3.1pt\lower-1pt
   \hbox{$\mathsurround=0pt\mathchar"0012$}\kern-4pt}\:$}}}
\newcommand{\drest}{\rest\rest}

\newcommand{\rs}{\!\!\!\!\!}
\newcommand{\pr}{\prec \!\!\!)}
\newcommand{\suc}{\succ \!\!(}
\newcommand{\nin}{\in\!\!\!\!\!/}
\newcommand{\rang}{<\!\!\! )}
\newcommand{\lang}{>\!\!\!\!\! (}

\newcommand{\empt}{\emptyset}

\newcommand{\cntd}{\subseteq}
\newcommand{\cnts}{\supseteq}
\newcommand{\pcntda}{\lower5pt\hbox{$\stackrel{\subset}{\neq}$}}
\newcommand{\pcntdb}{\lower5pt\hbox{$\stackrel{\supset}{\neq}$}}

\newcommand{\pexp}[1]{\hbox{$ #1 $}}

\newcommand{\real}[1]{\hbox{\rm #1}}
\newcommand{\Cal}[1]{{\cal #1}}

\newcommand{\ml}{\ell}

\newcommand{\ga}{\alpha}
\newcommand{\gga}{\gamma}
\newcommand{\gb}{\beta}
\newcommand{\gd}{\delta}
\newcommand{\gk}{\kappa}
\newcommand{\get}{\eta}
\newcommand{\gep}{\varepsilon}
\newcommand{\gvp}{\varphi}
\newcommand{\gl}{\lambda}
\newcommand{\gL}{\Lambda}
\newcommand{\Gl}{\Lambda}
\newcommand{\gch}{\chi}
\newcommand{\gp}{\pi}
\newcommand{\gps}{\psi}
\newcommand{\gs}{\sigma}
\newcommand{\gr}{\rho}
\newcommand{\ges}{\varsigma}
\newcommand{\gS}{\Sigma}
\newcommand{\gt}{\theta}
\newcommand{\gT}{\Theta}
\newcommand{\go}{\omega}
\newcommand{\gO}{\Omega}
\newcommand{\gG}{\Gamma}
\newcommand{\gx}{\xi}
\newcommand{\gz}{\zeta}
\newcommand{\e}{{\rm e}}
\newcommand{\m}{{\rm m}}

\newcommand{\bA}{\hbox{\bf A}}
\newcommand{\bC}{\hbox{\bf C}}
\newcommand{\bL}{\hbox{\bf L}}
\newcommand{\bI}{\hbox{\bf I}}
\newcommand{\bR}{\hbox{\bf R}}
\newcommand{\bE}{\hbox{\bf E}}
\newcommand{\bB}{\hbox{\bf B}}
\newcommand{\bN}{\hbox{\bf N}}
\newcommand{\bZ}{\hbox{\bf Z}}
\newcommand{\bT}{\hbox{\bf T}}
\newcommand{\bQ}{\hbox{\bf Q}}
\newcommand{\bU}{\hbox{\bf U}}
\newcommand{\bu}{\hbox{\bf u}}
\newcommand{\bl}{\hbox{\bf l}}
\newcommand{\bt}{\hbox{\bf t}}
\newcommand{\bh}{\hbox{\bf h}}
\newcommand{\bH}{\hbox{\bf H}}
\newcommand{\bP}{\hbox{\bf P}}
\newcommand{\bS}{\hbox{\bf S}}
\newcommand{\bq}{\hbox{\bf q}}

\newcommand{\HC}{\hbox{\rm\bf H$^{\infty}+$c}} 

\newcommand{\CA}{{\cal A}}
\newcommand{\CB}{{\cal B}}
\newcommand{\CC}{{\cal C}}
\newcommand{\CD}{{\cal D}}
\newcommand{\CE}{{\cal E}}
\newcommand{\CF}{{\cal F}}
\newcommand{\CG}{{\cal G}}
\newcommand{\CH}{{\cal H}}
\newcommand{\CI}{{\cal I}}
\newcommand{\CJ}{{\cal J}}
\newcommand{\CK}{{\cal K}}
\newcommand{\CL}{{\cal L}}
\newcommand{\CM}{{\cal M}}
\newcommand{\CN}{{\cal N}}
\newcommand{\CO}{{\cal O}}
\newcommand{\CP}{{\cal P}}
\newcommand{\CQ}{{\cal Q}}
\newcommand{\CR}{{\cal R}}
\newcommand{\CS}{{\cal S}}
\newcommand{\CT}{{\cal T}}
\newcommand{\CU}{{\cal U}}
\newcommand{\CV}{{\cal V}}
\newcommand{\CW}{{\cal W}}
\newcommand{\CX}{{\cal X}}
\newcommand{\CY}{{\cal Y}}
\newcommand{\CZ}{{\cal Z}}

\newcommand{\biU}{\hbox{$\bf\it U$}}
\newcommand{\rmn}[1]{{\rm (#1)\,\,}}

\newcommand{\bara}{\bar{a}}
\newcommand{\barA}{\bar{A}}
\newcommand{\barf}{\bar{f}}
\newcommand{\barh}{\bar{h}}
\newcommand{\bark}{\bar{k}}
\newcommand{\barn}{\bar{n}}
\newcommand{\barz}{\bar{z}}
\newcommand{\barkk}{\bar{\bar{k}}}
\newcommand{\barQ}{\bar{Q}}
\newcommand{\barga}{\bar{\ga}}
\newcommand{\bargep}{\bar{\gep}}
\newcommand{\bargb}{\bar{\gb}}
\newcommand{\bargl}{\bar{\gl}}
\newcommand{\bargd}{\bar{\gd}}
\newcommand{\bargo}{\bar{\go}}
\newcommand{\hatD}{\hat{D}}
\newcommand{\hatK}{\hat{K}}
\newcommand{\hatP}{\hat{P}}
\newcommand{\hatS}{\hat{S}}
\newcommand{\hatT}{\hat{T}}
\newcommand{\hatV}{\hat{V}}
\newcommand{\hatv}{\hat{v}}

\newcommand{\nuo}{\nu_{0}}

\newcommand{\veca}{\vec{a}}
\newcommand{\vecb}{\vec{b}}
\newcommand{\vecc}{\vec{c}}
\newcommand{\vecd}{\vec{d}}
\newcommand{\vecf}{\vec{f}}
\newcommand{\vecg}{\vec{g}}
\newcommand{\vecm}{\vec{m}}
\newcommand{\vecS}{\vec{S}}
\newcommand{\vecu}{\vec{u}}
\newcommand{\vecx}{\vec{x}}
\newcommand{\vecy}{\vec{y}}
\newcommand{\vecz}{\vec{z}}
\newcommand{\vecgf}{\vec{\gf}}
\newcommand{\vecgt}{\vec{\gt}}
\newcommand{\bgtu}{\bigtriangleup}
\newcommand{\tr}{\triangle}
\newcommand{\bgtd}{\bigtriangledown}
\newcommand{\bsl}{\backslash}
\newcommand{\nequ}{\equiv\!\!\!\!\! /}

\newcommand{\tl}{\tilde}
\newcommand{\tlb}{\tilde{b}}
\newcommand{\tlc}{\tilde{c}}
\newcommand{\tly}{\tilde{y}}
\newcommand{\tlG}{\tilde{G}}
\newcommand{\tlgt}{\tilde{\tau}}
\newcommand{\tlgw}{\tilde{\gw}}
\newcommand{\tlgs}{\tilde{\gs}}
\newcommand{\tlgm}{\tilde{\mu}}
\newcommand{\tlF}{\tilde{F}}
\newcommand{\tlJ}{\tilde{J}}
\newcommand{\tlgf}{\tilde{\phi}}
\newcommand{\tlY}{\tilde{Y}}
\newcommand{\tlgl}{\tilde{\lambda}}

\newcommand{\Raro}{\Rightarrow}
\newcommand{\LRaro}{\Leftrightarrow}
\newcommand{\raro}{\rightarrow}
\newcommand{\laro}{\leftarrow}
\newcommand{\Llraro}{\Longleftarrow}
\newcommand{\Lrraro}{\Longrightarrow}
\newcommand{\LLRraro}{\Longleftrightarrow}
\newcommand{\lngl}{\langle}
\newcommand{\rngl}{\rangle}

\newcommand{\bcap}{\bigcap}
\newcommand{\bcup}{\bigcup}
\newcommand{\sub}{\subset}
\newcommand{\cd}{\cdot}
\newcommand{\itms}[1]{\item[[#1]]}
\newcommand{\lnri}{\lim_{n\raro\infty}}
\renewcommand{\i}{\infty}

\newcommand{\ZZ}{{\mathchoice {\hbox{$\sf\textstyle Z\kern-0.4em Z$}}
{\hbox{$\sf\textstyle Z\kern-0.4em Z$}}
{\hbox{$\sf\scriptstyle Z\kern-0.3em Z$}}
{\hbox{$\sf\scriptscriptstyle Z\kern-0.2em Z$}}}}

\newcommand{\AAA}{{\mathchoice {\hbox{$\sf\textstyle A\kern-0.4em A$}}
{\hbox{$\sf\textstyle A\kern-0.4em A$}}
{\hbox{$\sf\scriptstyle A\kern-0.3em A$}}
{\hbox{$\sf\scriptscriptstyle A\kern-0.2em A$}}}}

\newcommand{\RR}{{\rm I\!R}}
\newcommand{\DD}{{\rm I\!D}}
\newcommand{\EE}{{\rm I\!E}}
\newcommand{\MC}{{\rm I\!\!\!C}}

\newcommand{\NN}{{\rm I\!N}}
\newcommand{\BB}{{\rm I\!B}}

\newcommand{\Cc}{{\mathchoice {\setbox0=\hbox{$\displaystyle\rm C$}\hbox
{\hbox
to0pt{\kern0.4\wd0\vrule height0.9\ht0\hss}\box0}}
{\setbox0=\hbox{$\textstyle\rm C$}\hbox{\hbox
to0pt{\kern0.4\wd0\vrule height0.9\ht0\hss}\box0}}
{\setbox0=\hbox{$\scriptstyle\rm C$}\hbox{\hbox
to0pt{\kern0.4\wd0\vrule height0.9\ht0\hss}\box0}}
{\setbox0=\hbox{$\scriptscriptstyle\rm C$}\hbox{\hbox
to0pt{\kern0.4\wd0\vrule height0.9\ht0\hss}\box0}}}}

\newcommand{\QQ}{{\mathchoice {\setbox0=\hbox{$\displaystyle\rm
Q$}\hbox{\raise
0.15\ht0\hbox to0pt{\kern0.4\wd0\vrule height0.8\ht0\hss}\box0}}
{\setbox0=\hbox{$\textstyle\rm Q$}\hbox{\raise
0.15\ht0\hbox to0pt{\kern0.4\wd0\vrule height0.8\ht0\hss}\box0}}
{\setbox0=\hbox{$\scriptstyle\rm Q$}\hbox{\raise
0.15\ht0\hbox to0pt{\kern0.4\wd0\vrule height0.7\ht0\hss}\box0}}
{\setbox0=\hbox{$\scriptscriptstyle\rm Q$}\hbox{\raise
0.15\ht0\hbox to0pt{\kern0.4\wd0\vrule height0.7\ht0\hss}\box0}}}}

\newcommand{\lun}[1]{\bigcup\limits_{#1}}

\newcommand{\II}{{\bf I\kern -1pt I}}

\newcommand{\half}{\frac{1}{2}}
\newcommand{\singcol}[2]{\left[\begin{rray}{c}#1\\#2\end{array}\right]}
\newcommand{\singcolb}[2]{\left(\begin{array}{c}#1\\#2\end{array}\right)}
\newcommand{\doubcol}[4]{\left[\begin{array}{cc}#1&#2\\#3&#4\end{array}\right]}
\newcommand{\doubcolb}[4]{\left(\begin{array}{cc}#1&#2\\#3&#4\end{array}\right)}
\newcommand{\sing}[1]{\langle #1\rangle}

\newcommand{\susim}{\stackrel{succ}{sim}}
\newcommand{\eqqdf}{\stackrel{def}{\equiv}}
\newcommand{\eqdf}{\stackrel{def}{=}}
\newcommand{\inj}{\stackrel{1-1}{\longrightarrow}}
\newcommand{\surj}{\vbox{\hbox{$\longrightarrow $
                  \kern -22pt \hbox{\lower 2.5pt  \hbox{\tiny onto}}
                  \kern -16pt \hbox{\raise 5pt  \hbox{\tiny 1-1}}
                  \kern 3pt}}}
\newcommand{\uarrow}[2]{\vbox{\hbox{$\longrightarrow $
                  \kern -16pt \hbox{\raise 5pt  \hbox{\tiny $#1$}}
                  \kern 10pt }}}
\newcommand{\AEQ}{$\forall \exists \:$ \,}
\newcommand{\A}{\forall}
\newcommand{\E}{\exists}

	\newcommand{\AH}{A\# H}

\newcommand{\R}{\hbox{I\kern-.1500em \hbox{\sf R}}}
\newcommand{\Q}
   {\hbox{${\rm Q} \kern -7.5pt \raise 2pt \hbox{\tiny$|$}\kern 7.5pt$}}
\newcommand{\C}
   {\hbox{${\rm C} \kern -7.5pt \raise 2pt \hbox{\tiny$|$}\kern 7.5pt$}}

\newcommand{\qd}{\kern 5pt\vrule height8pt width0.5pt depth0pt}

\newcommand{\bd}{\begin{description}}
\newcommand{\ed}{\end{description}}

\newcommand{\bc}{\begin{center}}
\newcommand{\ec}{\end{center}}

\newcommand{\beq}{\begin{equation}}
\newcommand{\eeq}{\end{equation}}

\newcommand{\ben}{\begin{enumerate}}
\newcommand{\een}{\end{enumerate}}

\newcommand{\seqn}[2]{\langle#1,\ldots ,#2\rangle}
\newcommand{\sseqn}[2]{\langle#1\,|#2\rangle}
\newcommand{\ssseqn}[2]{\langle#1\,|#2\rangle}
\newcommand{\setn}[2]{$\{\,$#1$\:|\:$#2$\}$}
\newcommand{\setm}[2]{\{\ #1\:|\:#2\}}

\newcommand{\setp}[2]{\{\,#1\: : \:#2\}}

\newcommand{\fsetn}[2]{\{\,#1,\ldots ,#2\}}
\newcommand{\ssetn}[2]{\Llbc\,#1\,|#2\Lrbc}
\newcommand{\sssetn}[2]{\LLlbc\,#1\,|#2\LLrbc}
\newcommand{\pair}[2]{\langle#1 ,#2\rangle}
\newcommand{\trpl}[3]{\langle#1 ,#2 ,#3 \rangle}
\newcommand{\fnn}[3]{#1:#2 \raro #3}
\newcommand{\surjfn}[3]{#1:#2 \surj #3}
\newcommand{\pgn}[3]{\langle#1 ,#2 ; \, #3 \rangle}
\newcommand{\iso}[3]{\m{ #1 : #2 \cong #3 }}

\newcommand{\hsetn}[2]{\hbox{\{\,#1\,|#2\}}}

\newcommand{\azc}{\hbox{$\aleph_{0}$-categorical }}
\newcommand{\azcy}{\hbox{$\aleph_{0}$-categoricity }}
\newcommand{\meq}{\hbox{$M^{Eq}$}}
\newcommand{\feq}{\hbox{$f^{Eq}$}}
\newcommand{\ut}{U-tree }
\newcommand{\uts}{U-trees }



\def\brd{\vrule height 2pt width .5pt depth 1.5pt}
\newcommand{\ol}{\overline}
\newcommand{\sm}{\setminus}
\newcommand{\minus}{^{-1}}
\newcommand{\bbraces}[1]{\left\{ #1 \right\}}


\def\newtheorems{\newtheorem{theorem}{Theorem}[section]
		 \newtheorem{cor}[theorem]{Corollary}
		 \newtheorem{prop}[theorem]{Proposition}
		 \newtheorem{lemma}[theorem]{Lemma}
		 \newtheorem{defn}[theorem]{Definition}
		 \newtheorem{claim}[theorem]{Claim}
		 \newtheorem{sublemma}[theorem]{Sublemma}
		 \newtheorem{example}[theorem]{Example}
		 \newtheorem{remark}[theorem]{Remark}
		 \newtheorem{question}[theorem]{Question}
		 \newtheorem{conjecture}{Conjecture}[section]}
\def\rar{\mathop{\longrightarrow}\limits}
\def\max{\mathop{\rm max}\limits}

\newcommand{\diaga}{\begin{picture}(40,12)
\put(24,20.2){\line(2,-1){69}}%
                   \end{picture}}

\newcommand{\thediaga}{\raisebox{-.2ex}{\makebox[0pt]{\diaga}}}

\newcommand{\rboxa}[1]{\raisebox{1.5ex}[2ex]{\hspace{.12in} #1}}

\newtheorem{theorem}{Theorem}
\newcounter{examp}
\font\bigcaps=cmcsc10 scaled 1250 %
\newcommand{\ExaLabel}[1]{\label{eqn:#1}}
\newcommand{\ExaRef}[1]{\ref{eqn:#1}}
\newcommand{\eop}[0]{\begin{flushright} $\blacksquare$ \end{flushright}}
\newtheorem{proposition}{Proposition}
\newtheorem{corollary}{Corollary}
\newtheorem{lemma2}{Lemma}
\newtheorem{remark}{Remark}
\newtheorem{example}{Example}
\newtheorem{question}{Question}
\newtheorem{definition}{Definition}
\newtheorem{property}{Property}
\newtheorem{conjecture}{Conjecture}

\newcommand{\BlackBox}{
\rule{1.1ex}{1.2ex} }
\newenvironment{Proof}[1]{\ifx\ShowProofs\True \addvspace{1.0ex}
{\bf Proof:}#1~\BlackBox\fi}{}
\def\proofi{\futurelet\next\lookforbracket}
\def\lookforbracket{\ifx\next[\let\go\usespecialterm
\else\let\go\relax
\vskip1sp\noindent{\bf Proof:}\enskip\relax\ignorespaces\fi\go}
\def\usespecialterm[#1]{\vskip12pt
\noindent{\bf #1:}\enskip\relax\ignorespaces}

\baselineskip 25pt
\title{\fontfamily{cmss} \selectfont Nonstationary  iterative processes}
\vskip 15pt

\author{\setcounter{footnote}{0}Tamara Kogan, \thanks{Amit Educational Network, Beer-Sheva, Israel}$\;\;$ \setcounter{footnote}{6}Luba Sapir, \thanks{Departments of Mathematics and Computer Science, Ben-Gurion University, Beer-Sheva, Israel,
e-mail: lsapir@bgu.ac.il}$\;\;$\setcounter{footnote}{4}Amir Sapir, \thanks{Department of Computer Science, Sapir Academic College, Sha'ar HaNegev, Yehudah, Israel and The Center for Advanced Studies in Mathematics, Ben-Gurion University,
Beer-Sheva, Israel, e-mail: amirsa@cs.bgu.ac.il}$\;\;$\setcounter{footnote}{5}Ariel Sapir \thanks{Department of Computer Science, Ben-Gurion University, Beer-Sheva, Israel,
e-mail: arielsa@cs.bgu.ac.il}}
\maketitle

\fontfamily{cmss} \selectfont
\abstract{\large In this paper we present iterative methods of high efficiency by the criteria of J. F. Traub and A. M. Ostrowski.

We define {\it s-nonstationary iterative processes} and prove that, for any one-point iterative process without memory, such as, for example, Newton's, Halley's, Chebyshev's methods, there exists an s-nonstationary process of the same order, but of higher efficiency.

We supply constructions of these methods, obtain their properties and, for some of them, also their geometric interpretation. The algorithms we present can be transformed into computer programs in straight-forward manner.  The methods are demonstrated by numerical examples.
\\



%
 }

Keywords: {\it one-point iterative method with memory, Traub-Ostrowski index of computational efficiency, informational efficiency, order of convergence, Kung-Traub conjecture.}

\large
\baselineskip 18pt


\section{\textnormal{\fontfamily{cmss} \selectfont Introduction}}~\label{intro} $\hskip 5pt$
\vskip-24pt 
Iterative methods for solving a single non-linear equation of the form
\begin{equation*}~\label{fx}
f(x)=0,
\end{equation*}
were known for a long time ago.  Originally, most popular processes were of form
\begin{equation}~\label{F1}
x_{i+1}=\varphi(x_{i}),\qquad i=0,1,2,\;\ldots,
\end{equation}
where $x_{0}$ is an initial approximation to a simple real root $\alpha$ of (\ref{fx}). Various schemes of methods, with emphasis on higher order of convergence, were devised. Some of these methods required increased computational cost per iteration. During the $60's$ of the last century, fundamental research was done by A.~M. Ostrowski \cite{Ostrowski1960} and J.~F. Traub \cite{Traub1964}. In these works, the issue of effectiveness of iterative methods was studied, and criteria for effectiveness were defined (see Subsection~\ref{ef}).

Traub~\cite{Traub1964} classified iterative methods by the information they require.
If $x_{i+1}$, for any $i$, is determined by the new information at $x_{i}$ and reused information at $x_{i-1}, x_{i-2}, \ldots, x_{i-k} \, (k=1, 2, \ldots, i)$, i.e.
\begin{equation}~\label{Fn4}
x_{i+1}=\varphi(x_{i}; x_{i-1}, x_{i-2},\ldots,x_{i-k}),\qquad 1 \le k \le i,
\end{equation}
then $\varphi$ is called a {\it one-point iterative function with memory}.
If $\varphi$ for any $i$ is determined only by the new information at $x_{i}$, i.e.
$$x_{i+1} = \varphi(x_{i})\, ,$$
it will be called a {\it one-point iterative function without memory}.

Recently, in \cite{KoganSapir2017} it was proved that the most efficient, in the class of one-point methods without memory, are those of order $3$, such as Halley's method (\ref{Halley}) and Chebyshev's method (\ref{Cheby}), and so on. These methods are subject matter of current study too, such as of Kumari and Parida \cite{KuPa2019}, which analyzed the local convergence for Chebyshev's method in Banach spaces by using majorizing sequence.

During $1960-1970's$, there was an increased use of multi-point iterative methods (with and without memory). Traub~\cite{Traub1964} proved that the computational efficiency of these methods is higher than that of one-point methods.

Note that all the above-mentioned iterative methods are stationary. Recall, that an iterative process is called {\it stationary} if its function is the same at each iteration. If the function depends on the iteration's number, then the process is called {\it nonstationary} (see \cite{OrRhe1970}).
During the last years, several works discussing nonstationary iterative methods appeared. These methods have higher efficiency, which can be explained, intuitively, as follows: For the first few approximations, which may still be far away from $\alpha$, it is less reasonable to apply iterative methods of high computation cost, so that not to ``use a sledgehammer to crack a nut". For nonstationary methods, both the efficiency and the computational cost increase as we get closer to the root.

In \cite{KoganSapir2007} the authors present a nonstationary iterative method, based on some class of stationary methods (named the Fibonacci family), studied in \cite{KoganSapir2016}.  In this paper we provide a geometric interpretation of the nonstationary process of \cite{KoganSapir2007}.

The methodology of construction of nonstationary methods was given in several papers.
Nasr Al Din \cite{NasrAdin2011} suggests an interesting nonstationary process, based on Halley's method. Jain and Sethi \cite{JainSethi2018} and Jain, Chand and Sethi \cite{JainChandSethi2018} combined Aitkin type methods with other iterative processes, and obtained nonstationary methods of high efficiency.

In this paper we define {\it s-nonstationary iterative processes},
 explore their properties, and supply applications.

The paper is organized as follows:
Section~\ref{def} contains the notations, definitions and some previous results, pertaining to this paper. Section~\ref{mainres} provides our main results, where we construct a nonstationary process of highest efficiency for the class of one-point processes with memory. We also present two examples of nonstationary processes. Section~\ref{num} contains numerical example. Section~\ref{conc} contains a comparison of the suggested processes with one-point and multi-point processes, brief conclusions, and some practical recommendations.

\section{\textnormal{\bf \fontfamily{cmss} \selectfont Definitions, notations and previous results}}~\label{def} $\hskip 5pt$
\vskip -24pt
\vskip -24pt

\subsection{\textnormal{{\bf\fontfamily{cmss} \selectfont Order of convergence}}}~\label{order of convergence}
Let $e_n=x_n-\ga$ be the {\it error} of the $n$-th approximation to $\alpha$. If
\begin{equation}~\label{en}
|e_{n+1}|=c|e_n|^k+o(|e_n|^k),
\end{equation}
where $c>0$ and $k\geq 1$ are some constants, then the iterative process {\it converges with order} $p=k,$ and the number $c$ is called the {\it asymptotic error constant}.


Wall~\cite{Wall1956} defined the order of the process in a more general way, namely by
\begin{equation}~\label{order}
p=\lim_{n \to \infty} \frac{p_{n+1}}{p_n},
\end{equation}
where $p_n=-\log |e_n|,$ if the limit of the right-hand side exists.


From the practical point of view, for sufficiently large $n,$ $p_n$ provides the number of correct places of
decimals in $x_n$ (with the $\log$ taken in base $10$). Thus $p_{n+1}\approx p \cdot p_n$, which means that the number of correct places of decimals
in $x_{n+1}$ is about $p \cdot p_n$. Roughly speaking, the order $p$ can be interpreted as the factor by
which the accuracy has been multiplied at each successive iteration.

Various representatives of higher-order of convergence iterative methods can be found in (\cite{AmBuPl2004}, \cite{CoTo2011}, \cite{GrDi2006b}, \cite{GuHe2001}, \cite{KuPa2019}, \cite{Kogan1966}, \cite{Kogan1967}, \cite{McNameePan2013}, \cite{SolaimanHashim2018}, \cite{SolaimanKarimHashim2018}, \cite{PhiriMakinde2018}, \cite{KungTr74}).\\

\subsection{\textnormal{{\bf\fontfamily{cmss} \selectfont Efficiency of iterative processes}}}~\label{ef}
There is a variety of situations where the calculation of the function $f$ or its derivatives turns out to be of very high computational cost, such that other calculations involved in the iterative process are negligible. For these situations, in $1960,$ Ostrowski~\cite{Ostrowski1960} suggested the unit of measure of computational work,  so-called {\it Horner unit}, which is equal to computational work involved in
evaluating a function or any of its derivatives. Suppose we need $d_n$ Horner units for passing from $x_{n}$ to $x_{n+1}$.
Ostrowski's methodology takes
\begin{equation}~\label{ind}
I=\lim_{n \to \infty} {p}^{\frac{1}{d_n}},
\end{equation}
(if the limit exists)
as the measure of the efficiency of the process -- its {\it efficiency index}. For a stationary process, $d_n$ is a constant, i.e., $d_n=d$, and the efficiency index $I=p^{\frac{1}{d}}$.
Analogously  to $p$, which roughly provides the rate of growth
of the number of correct digits per iteration, $I$ provides the corresponding rate per consumption of a
single Horner unit.

Approximately at the same years, Traub~\cite{Traub1964} suggested the following two measures of efficiency for iterative processes (in \cite{Traub1964} they appear as {\rm EFF} and {\rm EFF*}, here we use $I_{1}$ and $I_{2}$, respectively):
\begin{description}
  \item [(1)] {\it Informational Efficiency}\\
 is the order $p$ divided by the information usage $d$:
  \begin{equation*}~\label{EFF}
I_{1}={\frac{p}{d}},
\end{equation*}
where the information usage $d$ is the number of new pieces of information required per iteration (or, equivalently, it is the number of Horner units required per iteration).
  \item [(2)] {\it Computational Efficiency}
     \begin{equation*}~\label{EFF*}
I_{2}=p^{\frac{1}{d}},
\end{equation*}
which takes into account the ``cost" of calculating different functions and their derivatives. For stationary processes, the efficiency index $I$ of Ostrowski is identical to the computational efficiency $I_{2}$.
\end{description}

Recently, in $2012$  J.M. McNamee and V.Y. Pan \cite{McNameePan2012} suggested an improved estimate of efficiency for iterative processes based on  $I_2$, called\\
 \begin{description}
\item [(3)] {\it Local Efficiency}\\
 \begin{equation*}~\label{LEFF}
I_3=\frac{\log_{10} p}{d},
\end{equation*}
which is proportional to the inverse of the computational work needed to obtain a desired accuracy of solution.\\

According to \cite{McNameePan2012}, if after $n$ steps, the output error is bounded by $10^{-p^n}$ and must stay below a desired accuracy $10^{-D}$, then  $n\approx~\frac{\log_{10} D}{\log_{10} p}$. The computational work equals $$W=n\cdot d\approx \log_{10} D \cdot \frac{d}{\log_{10} p}.$$
Note that the efficiency of iterative method is proportional to inverse work:
\begin{equation*}~\label{work}
\frac{1}{W}\approx \frac{1}{\log_{10} D}\cdot \frac{\log_{10} p}{d}.
\end{equation*}
 Since $D$ is problem dependent, only the factor $I_3=\frac{\log_{10} p}{d}$ is connected with efficiency of iterative method. (In \cite{McNameePan2012} the authors denote this factor as {\it Eff}.) \\

Clearly, $I_3=\log_{10} I_2,$ therefore  $I_2$ and $I_3$ increase or decrease in the same time.
\end{description}

%
%


\subsection{\textnormal{{\bf\fontfamily{cmss} \selectfont Newton's polynomial and its derivative}}}~\label{dividedDifference}
The construction of many iterative processes is based on Newton's divided difference formula. Denote by
$$f_{k,k-s}=f(x_k,x_{k-1},\ldots,x_{k-s})=\left \{
\begin{array}{ll}
{\displaystyle f(x_k),}& s=0,\\\\
{\displaystyle \frac{f_{k,k-s+1}-f_{k-1,k-s}}{x_k-x_{k-s}}},& 1\leq s\leq k,
\end{array} \right. $$
a {\it divided difference of order} $s$ between $x_k,x_{k-1},\ldots,x_{k-s}$. For example,
$f_{k,k}=f(x_k),\;\;f_{k,k-1}=\frac{f(x_k)-f(x_{k-1})}{x_k-x_{k-1}},\;\;f_{k,k-2}=\frac{f_{k,k-1}-f_{k-1,k-2}}{x_k-x_{k-2}}$ and so
on. Let $x_j,\;\;0\leq n$ be approximations to $\alpha$. Choosing the points in decreasing order of the indices, i.e., $x_n, x_{n-1},\dots,x_0,$ we obtain Newton's polynomial:
\begin{equation}~\label{nwpol}
P_n(x)=f_{n,n}+\sum_{i=1}^{n} f_{n,n-i}\cdot\prod_{j=n-i+1}^{n}(x-x_j),
\end{equation}
therefore
\begin{equation}~\label{nddf}
f(x)=P_n(x)+R_{n}(x),
\end{equation}
where ${\displaystyle R_n(x)=\frac{f^{(n+1)}(\xi(x))}{(n+1)!}\cdot\prod_{j=0}^{n}(x-x_j)}$.
Denote by $$w_n(x)=\prod_{j=0}^{n}(x-x_j).$$
Using this notation, (\ref{nwpol}) is equivalent to:
\begin{equation}~\label{nwpol1}
P_n(x)=f_{n,n}+\sum_{i=1}^{n} f_{n,n-i}\cdot \frac{w_n(x)}{w_{n-i}(x)}.
\end{equation}
Since $w_n(x)=(x-x_n)\cdot w_{n-1}(x),$ we obtain
$$w^\prime_n(x)=w_{n-1}(x)+(x-x_n)\cdot w^\prime_{n-1}(x),$$
and therefore
$$w^\prime_n(x_n)=w_{n-1}(x_n).$$
Deriving (\ref{nwpol1}) and evaluating at $x_n$ we obtain:
\begin{equation}~\label{d1nwpol}
P^\prime_n(x_n)=\sum_{i=1}^{n} f_{n,n-i}\cdot \frac{w_{n-1}(x_n)}{w_{n-i}(x_n)}.
\end{equation}
Hence, deriving (\ref{nddf}), we have
\begin{equation}~\label{d1nddf}
f^\prime(x_n)= P^\prime_n(x_n)+R^\prime_{n}(x_n),
\end{equation}
where $R^\prime_{n}(x_n)=\frac{w_{n-1}(x_n)}{(n+1)!}\cdot f^{(n+1)}(\xi(x_n)).$
If $f^{(n+1)}$ is bounded in convergence region of iterative process,  then
\begin{equation}~\label{d1anddf}
f^\prime(x_n) \approx P^\prime_n(x_n).
\end{equation}
In the sequel, it will be convenient to follow these approximations.

\subsection{\textnormal{{\bf\fontfamily{cmss} \selectfont Previous results}}}~\label{prev_res}
Traub~\cite{Traub1964} considered two classes of one-point iterative functions with memory, namely: {\it interpolatory functions} and {\it derivative estimated functions}.

Traub~\cite{Traub1964} proved that the order of a one-point iterative process with memory is a unique real positive root of the equation
\begin{equation}~\label{lemmafk}
F_{n}(t) \equiv t^{n+1} - s \Sigma_{j=0}^{n} t^{j} = 0 \, ,
\end{equation}
where $s$ is the number of derivatives of $f(x)$ (including $f(x)$ itself), used in the iterative process, and $n$ is the number of points at which the old information is being reused.

It $\,\,$ is $\,\,$ easy$\,\,$  to $\,\,$ check$\,\,$  that $\,\,$ $F_{n}(s) = -\Sigma_{j=0}^{n} s^{j} < 0$ and that $F_{n}(s+1)=1>0$, i.e.
\begin{equation}~\label{Fsp}
s < p < s+1 \, ,
\end{equation}
and the informational usage or number of Horner units is $d = s$.
Therefore, for the index of efficiency (or informational efficiency)
$$I_{1} = \frac{p}{d} = \frac{p}{s}$$
it holds, by $(\ref{Fsp})$, that
\begin{equation}~\label{range}
 1 < I_{1} < 1 + \frac{1}{s} \le 2 \,.
\end{equation}

\section{\textnormal{\bf \fontfamily{cmss} \selectfont Main results}}~\label{mainres}
This section is composed of two subsections:

In Subsection~\ref{sec31new} the authors introduce definition of $s$-nonstationary process and determine its order of convergence.

In Subsection~\ref{sec32new} we prove that any one-point process without memory, such as Newton, Halley, Chebyshev and so on, can be improved by composing an appropriate nonstationary process of the same order, but with higher effectivity indexes. In this subsection we also present the general method of construction of these nonstationary processes.

 We also discuss the nonstationary process for $s=1,$ given by the authors in~\cite{KoganSapir2007}, and its connection with Newton's method.  Then we depict its geometric interpretation. In addition, for $s=2,$ we obtain two new nonstationary processes of order $3$ with high efficiency indices.

\subsection{\textnormal{\bf \fontfamily{cmss} \selectfont $\bf s$-nonstationary process and its convergence order}}~\label{sec31new}
Let $x_0,x_1$ be two initial approximations to $\alpha$. Consider the following one-point nonstationary process with memory:
\begin{eqnarray}~\label{nonstat_upd2}
\begin{tabular}{l}
 $x_{2} = \varphi_{1}(x_1,x_0),$\\
 $x_{3} = \varphi_{2}(x_2;x_{1},x_{0}),$ \\
 $\ldots$\\
 $x_{k+1} = \varphi_{k}(x_k;x_{k-1},x_{k-2},...,x_{0}), \qquad \,k=1,2, \ldots$. \\
\end{tabular}
\end{eqnarray}
where the iterative function $\varphi_{k}$, depends on the first $s$ derivatives of the function $f(x)$ (including $f(x)$ itself). Such iterative methods we call  {\it $s$-nonstationary processes}.
We will prove later that the convergence order of a $s$-nonstationary process is $s+1$. Note that, for an arbitrary $k,$ if we already have $(k+1)$ initial approximation to $\alpha$, than the $s$-nonstationary process (\ref{nonstat_upd2}) can be started
 by $x_{k+1} = \varphi_{k}(x_k;x_{k-1},x_{k-2},...,x_{0})$.\\
For any fixed $k$ the process
\begin{eqnarray}~\label{nonstat_withder_upd2}
\begin{tabular}{l}
 $x_{k+i} = \varphi_{k}(x_{k+i-1};x_{k+i-2},...,x_{i-1}), \qquad \,i=1,2, \ldots$. \\
\end{tabular}
\end{eqnarray}
is a stationary one-point process with memory. Denote its order of convergence by $r_{k}$.

An example of a $s$-nonstationary process with $s=1$ is the iterative method presented in~\cite{KoganSapir2007}.
By fixing $k$ in the nonstationary process of~\cite{KoganSapir2007} we receive the following stationary process:\\
For $k=1$
\begin{eqnarray}
\nonumber x_{i+1} = \varphi_{1}(x_{i};x_{i-1})=x_{i}-\frac{f_i}{f_{i,i-1}}, \qquad \,i=1,2, \ldots. \end{eqnarray}
It is the well-known secant method with initial approximations $x_0,x_1$, which has an order of convergence $r_{1}=\frac{1+\sqrt{5}}{2}\approx 1.618$.\\
For $k=2$
\begin{eqnarray}
\nonumber x_{i+2} &=& \varphi_{2}(x_{i+1};x_{i},x_{i-1})\\
\nonumber &=&x_{i+1}-\frac{f_{i+1}}{f_{i+1,i}+f_{i+1,i-1}(x_{i+1}-x_i)},\qquad i=1,2, \ldots. \end{eqnarray}
It is a generalized secant method with initial approximations $x_0,x_1,x_2$, which has an order of convergence $r_{2}\approx 1.84$ (cf.~\cite{Kogan1966}).

Similarly, various examples of stationary process (\ref{nonstat_withder_upd2}) with convergence order $r_k$ can be constructed for an arbitrary $s$-nonstationary process (\ref{nonstat_upd2}). Lemma~\ref{lem2} illustrates the relation between convergence orders of such processes as follows:

%

\begin{lemma2}~\label{lem2}
For any $k$ ($k=1,2,\ldots$): $\,$ $r_k < r_{k+1}$
\end{lemma2}

\begin{proofi}
Traub~\cite{Traub1964} proved that $r_k$ is a unique value positive root of the$\,\,$  equation$\,\,$  (\ref{lemmafk})$\,\,$
i.e. $\,\,$ $F_k(r_k) =0$ $\,\,$ for $\,\,$ any $\,\,$ $k$,$\,\,$  where$\,\,$  $F_k(t) = t^{k+1} - s\sum_{j=0}^{k}t^{k}$. Moreover, by (\ref{Fsp}) for any $k$ it holds that:
\begin{equation}~\label{Fsprk}
s < r_k < s+1 \, .
\end{equation}
From (\ref{lemmafk}) we have:
$$ F_k(t) - F_{k+1} (t) = t^{k+1} \left(\left(s+1 \right) -t\right).$$

For $t=r_{k+1}$ from the last equation and (\ref{Fsprk}) we have $ F_k(r_{k+1})>0$.
Since $F_k(s)<0$,$\,\,$ $F_k(r_k)=0$ and $r_k$ is a unique positive root of (\ref{lemmafk}), we obtain that $s<r_k <r_{k+1}$.

\eop
\end{proofi}

\begin{theorem}
The order of convergence of the nonstationary one-point iterative process $(\ref{nonstat_upd2})$ is $\,\,\,s+1$.
\end{theorem}

\begin{proofi}
By Weierstrass's theorem, a monotone bounded sequence $\{r_{k}\}$ has limit $r$, and by the properties of limits, it holds that
\begin{equation}~\label{sequence}
 1 \le s < r_{k} < r \le s+1 .
\end{equation}
We shall show that $r$ is the order of convergence of $(\ref{nonstat_upd2})$.

Indeed, the order of convergence of the nonstationary process $(\ref{nonstat_upd2})$ is
\begin{equation}~\label{orderp}
p=\lim_{k \to \infty} \frac{p_{k+1}}{p_k}=\lim_{k \to \infty} \frac{\log |x_{k+1}-\ga|}{\log |x_{k}-\ga|}.
\end{equation}
Since for a sufficiently large fixed $k$
$$x_{k+1} = \varphi_{k}(x_{k};x_{k-1},...,x_{0}),$$
we have:
$$|x_{k+1}-\ga|\approx c_k \cdot |x_{k}-\ga|^{r_k},$$
where $c_k>0$ is bounded. Thus
$$p=\lim_{k \to \infty} \frac{\log (c_k \cdot |x_{k}-\ga|^{r_k})}{\log |x_{k}-\ga|}=\lim_{k \to \infty}r_k=r.$$
Now we show that $r = s + 1$.

From equation (\ref{lemmafk}) for $t=r_k$ we have

$$ F_k(r_k) = {r_k}^{k+1} - s\cdot \frac{{r_k}^{k+1}-1}{r_k-1}=0 $$
or
$$ {r_k}^{k+2} - (s+1) \cdot {r_k} ^ {k+1} +s = 0$$

From the last equation we observe that:
\begin{equation}~\label{theoremrk}
 r_k = s + 1 -\frac{s}{{r_k}^{k+1}} .
\end{equation}
Since for any $k$ it holds that
$$1<r_1\leq r_k <r$$
$$ s+1 - \frac{s}{{r_1}^{k+1}} \leq r_k < s+1-\frac{s}{r^{k+1}}$$
When $k\rightarrow \infty$ we shall have that
$$ r = \lim_{k\rightarrow \infty} r_k = s+1$$
\end{proofi}
\begin{corollary}~\label{nonstat}
The nonstationary process $(\ref{nonstat_upd2})$ is more effective than any of the stationary processes $(\ref{nonstat_withder_upd2})$.
\end{corollary}

Indeed,
$$I_{1,k} = I_{1}(\varphi_{k}) = \frac{r_{k}}{s} < \frac{r}{s} = \frac{s+1}{s} = 1 + \frac{1}{s}$$
$$I_{2,k} = I_{2}(\varphi_{k}) = \sqrt[s]{r_{k}} < r^{\frac{1}{s}} = (s+1)^{\frac{1}{s}} \, . \qquad $$
$$I_{3,k} = \log _{10} I_{2,k} = \frac{1}{s} \cdot \log _{10} \left( s+1\right) \Box$$

\subsection{\textnormal{\bf \fontfamily{cmss} \selectfont Construction of new nonstationary processes with memory}}~\label{sec32new}
In this subsection we provide new nonstationary processes, with high efficiency indexes.

\ul{Case $s = 1$}:
In \cite{KoganSapir2007} the authors construct a nonstationary iterative process:
\begin{equation}~\label{100}
x_{k+1} =
x_{k} - \frac{f_{k}}{G_k}, \qquad k=1,2,\dots,
\end{equation}
where {$G_{k}(x)=\displaystyle \sum_{i=1}^{k} f_{k,{k-i}} \prod_{j=1}^{i-1}(x_k-x_{k-j})$,} and
 $f_{k,{k-i}}$ is the divided differences of order $i$.\\

In $\,\,$ \cite{KoganSapir2007}$\,\,$ it$\,\,$ was$\,\,$ proved$\,\,$ that $\,\,$the $\,\,$order$\,\,$ of (\ref{100})$\,\,$ equals$\,\,$ $2$,$\,\,$ and $\,\,$that$\,\,$  $I_{1} = I_{2} = 2$.



Thus, the nonstationary algorithm can be formalized as:
\begin{algorithm}
\caption{\texttt{Method (\ref{100})}} \label{newtonDiff}
\begin{algorithmic}[1]
\STATE Choose $x_0, x_1$.
\STATE Compute $f_0=f(x_0)$.
\STATE $k = 1$.
\REPEAT
\STATE Compute $f_k = f_{k,k} = f(x_k)$.\\
$\hskip -24pt$
/* Construct divided differences of order $i$, where $i=1, 2, \ldots, k$, compute $G_k$ and $x_{k+1}$ */
\STATE $\displaystyle{f_{k,k-i} = \frac{f_{k,k-i+1} - f_{k-1,k-i}}{x_{k} - x_{k-i}}, \qquad\qquad i=1,2,\ldots,k}$
\STATE $\displaystyle{G_k = f_{k,k-1} + \sum_{i=2}^{k}f_{k,k-i} \prod_{j=1}^{i-1}(x_k-x_{k-j})}$
\STATE $\displaystyle{x_{k+1} = x_{k} - \frac{f_k}{G_k}}$
\STATE k := k+1
\UNTIL required accuracy achieved
\STATE {\bf return} $x_{k}$
\end{algorithmic}
\end{algorithm}


Next we provide the geometric interpretation of the algorithm.

Let $P_{k}(x)$ be an interpolation polynomial of order $k$, constructed based upon the points $(x_{i},f_{i}), i=0,1,\ldots,k$.


 Recall that $f_k = P_k \left( x_k \right)$ and $G_k \left( x_k \right) = {P^{'}}_k \left(x _k \right)$ (see Subsection~\ref{dividedDifference}). Hence, from ($~\ref{100}$)
we obtain:
\begin{equation}~\label{101}
 x_{k+1} = x_{k} - \frac{P_{k}(x_{k})}{P^\prime_{k}(x_{k})}.
\end{equation}

i.e. $x_{k+1}$ is an approximation which is obtained from $x_k$ by Newton's Method for the equation $P_k \left( x \right) = 0$ ($k=1,\,2,\ldots$).

This process is illustrated by Figure~1. Starting with the initial approximations $x_{0}, x_{1}$  we obtain $x_2$ by the Secant method. The next approximation $x_3$ is the $x$-intercept of the tangent line $TP_2(x)$ to the graph of the interpolation polynomial $P_2(x)$.



\begin{figure}[ht]~\label{figure5}.
\begin{center}
\includegraphics[width=10.8cm]{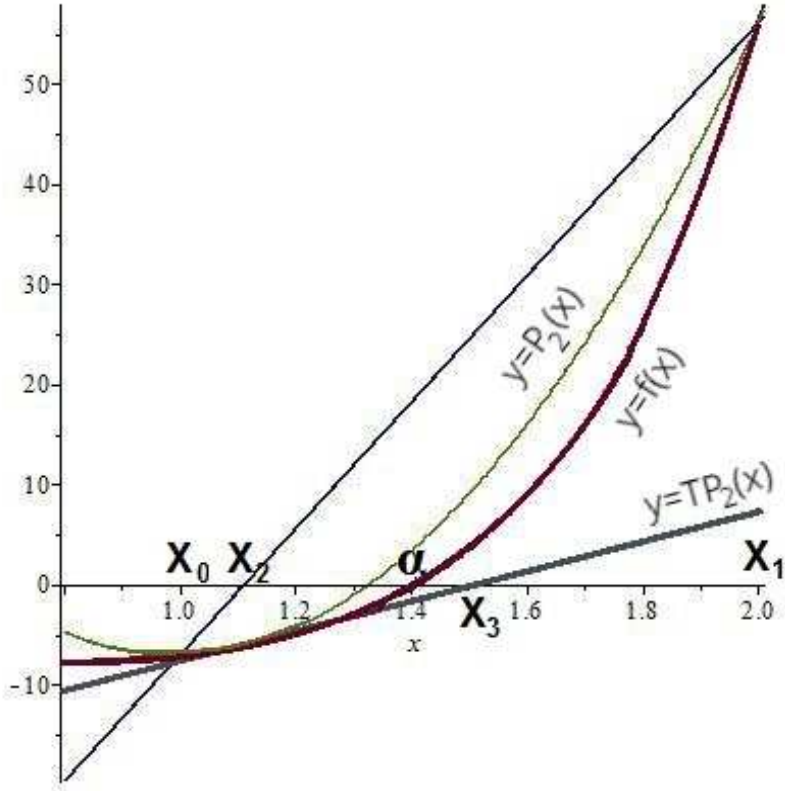}
\caption{Geometric interpretation of the nonstationary iterative process }
\end{center}
\end{figure}

Similarly, at the $i$-th step let $P_i(x)$ be the interpolation polynomial of the function $f(x)$ at the nodes $x_0,x_1,\ldots \,x_i$.  Thus, the next approximation $x_{i+1}$ to $\ga$ by the nonstationary iterative process is the $x$-intercept of the tangent line to the graph of the interpolation polynomial $P_i(x)$ at the point $(x_i,f(x_i))$.


\ul{Case $s = 2$}:

Next, we give a recipe for constructing a new nonstationary iterative method.

Consider an arbitrary iteration process with order of convergence $3$ without memory. For example, let us look at
Halley's method~\cite{Halley1964}:
\begin{equation}~\label{Halley}  x_{k+1} = x_{k} - \frac{2f_{k}f^{\prime}_{k}}{2(f^{\prime}_{k})^{2} - f_{k}f^{\prime\prime}_{k}},
\end{equation}

Let $x_0,\,x_1,\,x_2$ be initial approximations to $\alpha$. Denote $f^{\prime}(x) \equiv g(x)$; then $f^{\prime\prime}(x) \equiv g^{\prime}(x)$. For $k=1, 2, \ldots$, using points $(x_{i}, g_{i})$, where $g_{i} = g(x_{i}), i=0, 1, \ldots, k$, we construct the derivative of the interpolation polynomial of order $k$:
\begin{equation}~\label{DeriveIntPol}
 G_k=P^{\prime}_{k}(x_k) = w_{k-1}(x_k)\cdot \sum_{i=1}^{k}\frac{g_{k,k-i}}{w_{k-i}(x_k)}, \,\,\,\,\,k=2,\,3,\ldots
\end{equation}
where $g_{k,k-i}$ are the divided differences of order $i$.
Then we substitute $f_{k}^{\prime\prime}$ in ($\ref{Halley}$) by $P^{\prime}_{k}(x)$, obtaining the following nonstationary iterative process
\begin{equation}~\label{HalleyNonStatItProc01}  x_{k+1} = x_{k} - \frac{2f_{k}g_{k}}{2 g_k^{2} - f_{k}G_k},
\end{equation}
For $k \ge 2$, at any step, only two function evaluations are required: $f_{k}$ and $g_{k}$. Thus, $d = s = 2$ and the order of ($\ref{HalleyNonStatItProc01}$) is $p = s+1 = 3$.

A natural question arises: Is it beneficial to replace Halley's method by the nonstationary process $(\ref{HalleyNonStatItProc01})$ of the same order? In order to answer that, we compute the efficiency indexes of each one of them. For Halley's method:
$$p=3, \,\,\,s=d=3, \,\,\,I_{1}=\frac{p}{d}=1, \,\,\,I_{2}=p^{\frac{1}{d}}=3^{\frac{1}{3}},  \,\,\,I_{3}= \frac{\log _{10} 3}{3},$$
whereas for the method $(\ref{HalleyNonStatItProc01})$:
$$p=3,\,\,\, s=d=2, \,\,\,I_{1}=\frac{p}{d}=1.5, \,\,\,I_{2}=p^{\frac{1}{d}}=3^{\frac{1}{2}}, \,\,\,I_{3}= \frac{\log _{10} 3}{2} .$$

The nonstationary method $(\ref{HalleyNonStatItProc01})$ can be formalized as:
\begin{algorithm}
\caption{\texttt{Procedure Non-stationary Halley}} \label{halleyDiff}
\begin{algorithmic}[1]
\STATE Choose $x_0, x_1, x_2$.
\STATE Compute $f_0=f(x_0), f_1=f(x_1), g_0=g(x_0), g_1=g(x_1)$. /* $g(x) = f'(x)$ */
\STATE $k = 2$.
\REPEAT
\STATE Compute $f_k=f(x_k), g_k=g(x_k)$.
\STATE Construct divided differences of $g(x)$, compute $G_k$ and $x_k$
\STATE $\displaystyle{g_{k,k-i} = \frac{g_{k,k-i+1} - g_{k-1,k-i}}{x_{k} - x_{k-i}}, \qquad\qquad i=1,2,\ldots,k}$
\STATE $\displaystyle{G_k = w_{k-1}(x_k)\cdot \sum_{i=1}^{k}\frac{g_{k,k-i}}{w_{k-i}(x_k)}, \;\;\;\; {\rm where \;\;}  w_m(x)=\prod_{j=0}^{m}(x-x_j)}$
\STATE $\displaystyle{x_{k+1} = x_{k} - \frac{2 f_k g_k}{2g_k^2 - f_k G_k}}$
\STATE k := k+1
\UNTIL required accuracy achieved
\STATE {\bf return} $x_{k}$
\end{algorithmic}
\end{algorithm}

Similarly, one can take Chebyshev's method, of order 3:
\begin{equation}~\label{Cheby}
x_{k+1} = x_{k} - \frac{f(x_{k})}{f^{\prime}(x_{k})} \cdot \left(1 + \frac{f(x_{k})f^{\prime\prime}(x_{k})}{2\left( f^{\prime}(x_{k}) \right)^{2}}\right),
\end{equation}
and construct the nonstationary process
\begin{equation}~\label{ChebyNonStatItProc02}
x_{k+1} = x_{k} - \frac{f(x_{k})}{f^{\prime}(x_{k})} \cdot \left(1 + \frac{f(x_{k})G_k)}{2\left( f^{\prime}(x_{k}) \right)^{2}}\right).
\end{equation}
The effectivity indices of (\ref{Cheby}) and (\ref{ChebyNonStatItProc02}) are, respectively, those of (\ref{Halley}) and (\ref{HalleyNonStatItProc01}).

\ul{Case $s \ge 3$}:

Let $x_{k+1} = F(x_{k})$ be any one-point process without memory of order $p = s+1$.
Let the iterative function $F(x)$ depend on $f(x), f^{'}(x), \ldots, f^{(s)}(x)$. Similarly to case $s\!=\!2$, denote $f^{(s\!-\!1)}(x) \!\!\equiv\!\! g(x)$; then $f^{(\!s\!)}(x) \!\!\equiv\!\! g^{\prime}(x)$.
For $k=1, 2, \ldots$, using the points $(x_{i}, g_{i})$, where $g_{i} = g(x_{i}), i=0, 1, \ldots, k$, we construct the derivative of the interpolation polynomial of order $k$ and substitute $f^{(s)}(x_{k})$ in $F(x_{k})$ by $P^{\prime}_{k}(x_k)$, obtaining a nonstationary process of order $s+1$.

The efficiency indices of the process $x_{k+1} = F(x_{k})$ are
$$I_{1}=\frac{p}{d}=\frac{s+1}{s+1}=1, \,\,\,\,I_{2}=p^{\frac{1}{d}}=(s+1)^{\frac{1}{s+1}}, \,\,\,\,I_{3}=\frac{1}{s+1} \cdot \log _{10} \left(s+1\right),$$
whereas the nonstationary process' indices are:
$$I_{1}=\frac{p}{d}=\frac{s+1}{s}=1+\frac{1}{s}, \,\,\,\,I_{2}=p^{\frac{1}{d}}=(s+1)^{\frac{1}{s}}, \,\,\,\, I_{3}=\frac{1}{s} \cdot \log _{10} \left(s+1\right).$$
Hence, from the above and by Corollary~\ref{nonstat} we obtain:
\begin{theorem}~\label{thm2}
For any one-point iterative process, one can construct a more efficient $s$-nonstationary one-point iterative process of the same order.
\end{theorem}

\section{Numerical Results}~\label{num}

The following example illustrates the suggested methods $(\ref{HalleyNonStatItProc01})$ and $(\ref{ChebyNonStatItProc02})$. Consider the equation\\
\begin{equation}~\label{dugma}
f(x)=x^2-e^{\frac{1}{x}\cdot{\sin{\frac{\pi\cdot x^2}{2}}}}-1=0,
\end{equation}
which has $\ga=\sqrt{2}$ as a simple root. The following table illustrates the computation by formula $(\ref{HalleyNonStatItProc01})$ and $(\ref{ChebyNonStatItProc02})$, respectively, starting with $x_0=1.7,\;\;x_1=1.6,$ and $x_2=1.5$.
The correct value of the root $\ga$ to $10$ decimal places is $1.4142135624$.\\
\newcommand{\rbox}[1]{\raisebox{1ex}[2.3ex]{\hspace{.005in} #1}}

\begin{table}[ht]~\label{tb0}
  \begin{center}
  \begin{tabular}{|l|l|l|l|l|}
    \hline
    \multirow{2}{*}{$i$} &
      \multicolumn{2}{c|}{Method \;(\ref{HalleyNonStatItProc01})} &
      \multicolumn{2}{c|}{Method \;(\ref{ChebyNonStatItProc02})} \\
    & \rbox{$\!\!\!x_i\!\!\!$} &\rbox{$\!\!\!|e_i|=|x_i - \ga|\!\!\!$} & \rbox{$\!\!\!x_i\!\!\!$} &\rbox{$\!\!\!|e_i|=|x_i - \ga|\!\!\!$} \\
    \hline
    $0$ & $1.7$ & $0.2857864376$ & $1.7$ & $0.2857864376$ \\
    \hline
    $1$ & $1.6$ & $0.1857864376$  &$1.6$ & $0.1857864376$  \\
    \hline
   $2$ & $1.5$ & $0.0857864376$  &$1.5$ & $0.0857864376$ \\
    \hline
    $3$ & $1.4143581722$ & $0.0001446099$  & $1.4149666839$ & $0.0007531215$ \\
    \hline
    $4$ & $1.4142135632$ & $0.0000000008$  & $1.4142135854$ & $0.0000000009$ \\
    \hline
    $5$ & $1.4142135623$ & $2.98\cdot 10^{-62}$  & $1.4142135623$ & $2.02\cdot 10^{-62}$ \\
    \hline
    $6$ & $1.4142135624$ & $0$  & $1.4142135624$ & $0$ \\
    \hline
  \end{tabular}
  \end{center}
  \caption{Illustration of the suggested nonstationary methods.}
\end{table}

\section{\textnormal{\fontfamily{cmss} \selectfont Summary}}~\label{conc}
\vskip -24pt
The main conclusions of the paper are:
\begin{enumerate}
\item For $s$-nonstationary methods, the order of convergence is $p = s+1$ and the efficiency indices are
 $$I_{1} = 1 + \frac{1}{s}, \qquad I_{2} = (s+1)^{\frac{1}{s}},\qquad I_{3} = {\frac{\log_{10}(s+1)}{s}}.$$
 Obviously, the most effective among all $s$-nonstationary processes are those of $s = 1$, for which
 $$p = 2, \qquad I_{1} = I_{2} = 2\;\; {\rm and}\;\; I_{3} = \log_{10}{2}=0.301.$$

\item Theorem~\ref{thm2} yields that for any one-point iterative process (with or without memory) exists
 a $s$-nonstationary process of the same order, but more effective.

\item According to Kung-Traub conjecture the most effective methods among all multi-point iterative processes without memory  are methods of order $2^{n-1},$ where $n$ is number of function evaluations per iteration. Computational efficiency of these methods is $I_2=2^{\frac{n-1}{n}}<2$ for any $n$. Hence, a $s$-nonstationary iterative method for $s=1$ is more effective than any multi-point method without memory.

\item Efficiency index is a conditional term, since it is measured by Horner units. In practice, the amount of computations involved in evaluation may vary considerably from function to function. Thus, the nonstationary methods we offer are effective in particular if the evaluations of the derivatives require a large amount of computations, or are non-existent, or if the function is given by a table of values.

%

 \end{enumerate}

\bibliographystyle{elsart-num-sort}
\baselineskip 10pt
\normalsize

\end{document}